# Asymptotic of the greatest distance between adjacent primes and the Hardy-Littlewood conjecture

VICTOR VOLFSON

ABSTRACT   The paper substantiates the conjecture of the asymptotic behavior of the largest distance between consecutive primes: $\sup_{p_l \leq x}(p_{l+1} - p_l) \Box 2e^{-\gamma}\log^2(x)$, where $\gamma$ is the Euler constant. The Hardy-Littlewood conjecture about the number of prime tuples is investigated and the rationale for this conjecture is given taking into account the dependence of events that a large natural number is not divided into different prime numbers. It also substantiates why the accuracy of this conjecture is not affected by another assumption about the probability of a natural number being prime, although such a probability does not exist. The paper also considers the distribution of prime tuples using a mathematical model based on the Hardy-Littlewood conjecture.

1. INTRODUCTION

The opinion of mathematicians on probability theory was different at the beginning of the last century. This was due to the fact that at that time probability theory was not considered a branch of mathematics. Well-known mathematicians Hardy and Littlewood even wrote about this in their work – «Probability is not a notaition of pure mathematics, but of philosophy or physics».

However, despite this, the Hardy-Littlewood hypothesis about the number of prime tuples, based on a probabilistic approach, demonstrates high accuracy [1].

The high accuracy of the probabilistic model in the Hardy-Littlewood conjecture has already prompted William Banks, Kevin Ford, and TerenceTao to return to the probabilistic model of primes in a new article [2] to consider the greatest distance between adjacent primes. We will return to this model later. Now a little about the history of this issue.

The probabilistic model of primes was used by Kramer to estimate the greatest distance between adjacent primes in 1936. He proved in his work [3], that almost everywhere, under certain assumptions (see below), it holds: $\limsup_{x \to \infty} \frac{G_{kr}(x)}{\log^2 x} = 1$, where $G_{kr}(x)$ is the distance between adjacent elements of the set of the $kr$ Cramer model.

Kramer suggested in his probabilistic model, that each natural number can be included in a set $kr$ with probability $1/\log n$ and these events are mutually independent. One of the drawbacks of this model is that any natural number can be in one of the residue classes by prime modulo $p: 1, 2, ..., p-1$, the natural numbers from the set $kr$ are evenly distributed over all residue classes modularly $p$ in the Cramer model. Based on this probabilistic model, Kramer hypothesized that the greatest distance between adjacent primes is determined by the following asymptotic formula:

$$\sup_{p_l \leq x}(p_{l+1} - p_l) \Box \log^2 x, \qquad (1.2)$$

where $p_l, p_{l+1}$ are adjacent primes.

Grenville (in his probabilistic model of primes [4]) avoided the above drawback of the Cramer model. He put the parameters $A = \log^{1-Q(1)} x$ and $Q = \prod_{p \leq A} p$ for each interval $(x, 2x]$, where $x$ is the degree 2. Then he removed $n$ for which $(n, Q) > 1$ and he included $g$ in the set each remaining natural $n \in (x, 2x]$ with probability $\frac{Q}{\varphi(Q)\log(n)}$, where $\varphi$ is the Euler function. However, all events remains jointly independent in this model, as well as in the Cramer model. Based on his probabilistic model, Grenville makes the conjecture of the greatest distance between adjacent primes:

$$\sup_{p_l \leq x}(p_{l+1} - p_l) \geq (1 - o(1))\xi \log^2 x, \ (x \to \infty) \qquad (1.3)$$

where $\xi = 2e^{-\gamma} = 1,229..., \gamma$ is the Euler constant.

A new probabilistic model of prime numbers is considered in [2], which is based on the Hardy-Littlewood conjecture on the number of prime tuples mentioned above. The authors think that the new probabilistic model of primes:

$$R = \{n \geq e^2, n \in S_z\}, \qquad (1.4)$$

where a random set $S_z = z \setminus \bigcup_{p \leq z}(a_p \mod p)$ will be useful for studying the greatest distance between primes. Based on this model, the authors propose a conjecture about the asymptotic behavior of the largest distance between primes:

$$g((\xi - o(1))\log^2 x) \leq \sup_{p_l \leq x}(p_{l+1} - p_l) \leq g((\xi + o(1))\log^2 x), \ (x \to \infty) \qquad (1.5)$$

where $\xi = 2e^{-\gamma} = 1,229..., \gamma$ is the Euler constant.

Aauthors write that it seems to them that $g(a) \square g(b)$ when $a \square b$ in relation to the function $g$, although they cannot prove it. Authors suggest that conjecture (1.5) can be written in a more compact form:

$$\sup_{x \leq p_l}(p_{l+1} - p_l) \square g(\xi \log^2 x) \qquad (1.6)$$

and consider that in (1.6) $g(u) \square u$ with $u \to \infty$ and the conjecture can be written in the form:

$$\sup_{x \leq p_l}(p_{l+1} - p_l) \square \xi \log^2 x . \qquad (1.7)$$

We will try to substantiate the conjecture (1.7) in the second chapter of the paper.

The Hardy-Littlewood conjecture about the number of prime tuples (discussed above) is considered in [5]. This conjecture assumes independence in the totality of events that a large natural number is not divided into prime numbers, although this is not true, as shown in this paper. The proof of this fact and the substantiation of the Hardy-Littlewood conjecture of prime tuples taking this fact into account are given in the third chapter of the work. It also provides a rationale for why the conjecture accuracy is not affected by another assumption about the probability of a natural number being prime, although such a probability does not exist.

The fourth chapter of the paper is devoted to the search for the distribution of prime tuples using a mathematical model based on the Hardy-Littlewood conjecture.

## 2. ASYMTOTICS OF THE GREAT DISTANCE BETWEEN CONSECTIVE PRIMES

Let $f$ be a sequence of primes. Denote by $P_{a,b}$ the discrete uniform measure on the interval $[a,b)$. It is known that this measure is a probabilistic measure and for a sequence of primes is equal to the density of the sequence $f$ on this interval - $d(f,a,b)$.

Lemma 2.1

The probability of a natural number from the interval $[2, x)$ to be a prime is equal to:

$$P_{2,x}(f) = d(f,2,x) = 1/\log(x)(1+o(1)). \tag{2.1}$$

Proof

Based on the asymptotic law of primes, the number of primes not exceeding $x$ is equal to:

$$\pi(x) = x/\log(x)(1+o(1)) \tag{2.2}$$

at $x \to \infty$.

Having in mind (2.2) we obtain the value of the density of primes in the interval $[2, x)$:

$$d(f,2,x) = \pi(x)/x = 1/\log(x)(1+o(1)) \tag{2.3}$$

Based on the fact that the density found in (2.3) is a probability measure, we obtain:

$$P_{2,x}(f) = d(f,2,x) = 1/\log(x)(1+o(1)),$$

which corresponds to (2.1).

Theorem 2.2

The following asymptotic formula holds for the greatest distance between consecutive primes:

$$\sup_{p_l \leq x}(p_{l+1} - p_l) \square\ 2e^{-\gamma} \log^2(x), \tag{2.4}$$

where $\gamma$ is the Euler constant.

Proof

Following Grenville [4] the number of primes in the interval is determined using the sieve of Eratosthenes according to the formula:

$$\prod_{p \leq \sqrt{x}} (1 - 1/p) x \text{ with an error } R(x).$$

Thus, the density of a sequence of primes in the interval $[2, x)$ is determined by the formula:

$$d(f, 2, x) = \frac{\prod_{p \leq \sqrt{x}} (1 - 1/p) x + R(x)}{x} = \prod_{p \leq \sqrt{x}} (1 - 1/p) + \frac{R(x)}{x}. \tag{2.5}$$

Based on Lemma (2.1), we find the probability of a natural number from the interval $[2, x)$ to be a prime on the right in formula (2.5).

The following variants of residuals from division a natural number $x$ by a prime number $p$ are possible: $0, 1, ..., p-1$. Moreover, if the natural number $x$ is not divisible by $p$, there will only $p-1$ variants: $1, 2, ..., p-1$. Based on the Dirichlet theorem, in each class of residues from division $x$ by a prime number $p$ there will be the same number of natural numbers (if value $x \to \infty$), therefore the probability of the event that a natural number $x$ is not divided by a prime number $p$ if value $x \to \infty$ is equal $\frac{p-1}{p} = 1 - 1/p$.

We denote the event $A$ that a positive integer $x$ is not divisible by primes $p \leq \sqrt{x}$. Having in mind that the probability of an event that a natural number $x$ is not divisible by a prime number $p$ is $1 - 1/p$ (when $x \to \infty$). Assuming (as in the Cramer conjecture) the independence of these events, we get that the probability of the event $A$ is:

$$\Pr(A) = \prod_{p \leq \sqrt{x}} (1 - 1/p). \tag{2.6}$$

On the other hand, based on Lemma 2.1, the value - $\Pr(A) = P_{2,x}(f)$, therefore, based on (2.6), we obtain:

$$1/\log(x)(1 + o(1)) = \prod_{p \leq \sqrt{x}} (1 - 1/p) \tag{2.7}$$

or $1/\log(x) \Box \prod_{p \leq \sqrt{x}} (1-1/p)$.

Based on (2.5), events that a natural number $x$ is not divisible by different primes $p \leq \sqrt{x}$ are not independent events. This is confirmed by the Mertens theorem:

$$\prod_{p \leq \sqrt{x}} (1-1/p) \Box \frac{2e^{-\gamma}}{\log(x)}, \tag{2.8}$$

where $\gamma$ is the Euler constant.

Having in mind and according to the Cramer hypothesis, the maximum distance between adjacent primes is determined by the asymptotic formula:

$$\sup_{p_l \leq x} (p_{l+1} - p_l) \Box \log^2 x,$$

then based on (2.7) and (2.8) we obtain (2.4).

Thus, we took into account the dependence of events that a natural number is not divided into different primes $p \leq \sqrt{x}$, to obtain (2.4) in the probabilistic model of primes.

Theorem 2.2 is a justification of conjecture (1.7).

3. HARDY-LITTLEWOOD CONJECTURE ABOUT THE NUMBER OF PRIME TUPLES

Lemma 3.1

Let $x$ is a natural number, then events that $x$ cannot be divided into primes numbers $p \leq \sqrt{x}$ are dependent with a dependence coefficient $C = 0,5e^{\gamma}$ (when the value $x \to \infty$), where $\gamma$ is the Euler constant.

Proof

Having in mind (2.5) we have:

$$d(f, 2, x) = \prod_{p \leq \sqrt{x}} (1-1/p) + \frac{R(x)}{x} = C(x) \prod_{p \leq \sqrt{x}} (1-1/p), \tag{3.1}$$

where $C(x) = 1 + \dfrac{R(x)}{x \prod_{p \leq \sqrt{x}} (1 - 1/p)}$.

Based on Lemma 2.1, expression (3.1) is a probability $P_{2,x}(f)$ and it can be written in the form:

$$P_{2,x}(f) = C(x) \prod_{p \leq \sqrt{x}} (1 - 1/p) = 1/\log(x)(1 + o(1)) \qquad (3.2)$$

or $1/\log(x) \square C \prod_{p \leq \sqrt{x}} (1 - 1/p)$, where $C = \lim_{x \to \infty} C(x)$.

Based on the Mertens theorem, the value $C = 0,5e^{\gamma}$ in expression (3.2).

Now we justify the Hardy-Littlewood conjecture for prime twins, taking into account Lemma 3.1. We will use the assumption of the Hardy-Littlewood conjecture that the probability of an event $A_1$ that a large positive integer $x$ is prime exists and is equal to:

$$\Pr(A_1) = 1/\log(x). \qquad (3.3)$$

The authors of the conjecture do not consider that (3.3) is true, but simply assume.

Based on Lemma 2.1, we can only say that the probability of a natural number in the interval $[2, x)$ to be prime is equal to:

$$P_{2,x}(f) = 1/\log(x)(1 + o(1)).$$

Taking into account the premise of conjecture (3.3) we denote $A_2$ - the event that a natural number $x + 2$ is prime.

Now we formulate the Hardy-Littlewood conjecture for prime twins based on the above notation.

Conjecture 3.2

$$\Pr(A_1 \square A_2) \square C_2 / \log^2 x, \qquad (3.4)$$

where $C_2 = 2 \prod_{p \leq \sqrt{x}} \dfrac{p(p-2)}{(p-1)^2}$.

Proof

Having in mind (3.3) and the Mertens theorem, we can write:

$$\Pr(A_1) = 1/\log(x) \approx 0{,}5 e^{\gamma} \prod_{p \geq \sqrt{x}} (1 - 1/p). \tag{3.5}$$

Based on Lemma 3.2 and taking into account (3.3), we obtain:

$$\Pr(A_2) = 1/\log(x+2) \approx 0{,}5 e^{\gamma} \prod_{p \leq \sqrt{x+2}} (1 - 1/p). \tag{3.6}$$

An increase in the number 2 does not affect the asymptotic; therefore, (3.6) can be written as:

$$\Pr(A_2) = 1/\log(x) \approx 0{,}5 e^{\gamma} \prod_{p \geq \sqrt{x}} (1 - 1/p). \tag{3.7}$$

Remaining within the framework of (3.3), the probability that the natural numbers $x, x+2$ are simultaneously prime, having in mind (3.5), (3.7) and the dependence of these events, we obtain:

$$\Pr(A_1 \cap A_2) = \Pr(A_1)\Pr(A_2/A_1) = \Pr(A_1) C_2(x) \Pr(A_2) = C_2(x)/\log^2 x, \tag{3.8}$$

where $C_2(x)$ is the coefficient taking into account the dependence of these events.

Now we find the coefficient $C_2(x) = \dfrac{\Pr(A_2/A_1)}{\Pr(A_2)}$.

We divide all natural numbers (for this purpose) into classes modulo a prime number $p$. We obtain $p$ classes of natural numbers in this case, each of which will have the same number of numbers. It must be one of $p-1$ classes for a natural number which is not multiple $p$. There are only $p-2$ classes of these, where the natural number $p+2$ is not multiple $p$. Therefore, the probability of an event that natural numbers $x, x+2$ are not multiple $p$ is equal to:

$$\frac{p-2}{p-1}. \tag{3.9}$$

Having in mind (3.9) and the dependence of events that a natural number $p+2$ is not divisible by different primes (based on Lemma 3.1), we obtain:

$$\Pr(A_2 / A_1) \Box \, 0{,}5e^{\gamma} \prod_{p \leq \sqrt{x}} \frac{p-2}{p-1}. \tag{3.10}$$

Based on formulas (3.5), (3.7) and (3.10) we obtain:

$$C_2(x) = \frac{\Pr(A_2 / A_1)}{\Pr(A_2)} \Box \frac{0{,}5e^{\gamma} \prod_{p \leq \sqrt{x}} \frac{p-2}{p-1}}{0{,}5e^{\gamma} \prod_{p \leq \sqrt{x}} 1-1/p} = 2 \prod_{p \leq \sqrt{x}} \frac{p(p-2)}{(p-1)^2}, \tag{3.11}$$

which corresponds to (3.4).

Thus, the coefficients $0{,}5e^{\gamma}$ in formula (3.11) are reduced and do not affect the final result.

Now we formulate the Hardy-Littlewood conjecture in a general form. Moreover, assumption (3.3) remains valid.

Conjecture 3.3

$$\Pr(A_1 \Box .. \Box A_m) \Box \, C_m / \log^m x, \tag{3.12}$$

where $A_l$ is the event that the positive integer $x + 2k_1 + ... + 2k_l$ ($k_1 < k_2 < ... < k_l$) are primes, and

$$C_m = \frac{\prod_{p \leq x} 1 - \omega_m / p}{\prod_{p \leq x} (1 - 1/p)^m}.$$

Proof

Having in mind (3.3) and the Mertens theorem, we can write:

$$\Pr(A_1) = 1/\log(x) \Box \, 0{,}5e^{\gamma} \prod_{p \geq \sqrt{x}} (1 - 1/p). \tag{3.13}$$

Based on Lemma 3.2 and taking into account (3.3), we obtain:

$$\Pr(A_l) = 1/\log(x + 2k_1 + ... + 2k_l) \Box \, 0{,}5e^{\gamma} \prod_{p \leq \sqrt{x + 2k_1 + ... + 2k_l}} (1 - 1/p), \tag{3.14}$$

where $l = 1, ..., m-1$.

Since the increase in the number $2k_1+...+2k_l$ does not affect the asymptotic, then we can write (3.14) in the form:

$$\Pr(A_l) = 1/\log(x) \square\, 0,5 e^{\gamma} \prod_{p \geq \sqrt{x}} (1-1/p). \qquad (3.15)$$

Remaining within the framework of (3.3), the probability that the positive integers $x, x+2k_1,..., x+2k_1+...+2k_l$ are simultaneously primes, having in mind (3.13), (3.15) and the dependence of these events, we obtain:

$$\Pr(A_1 \square .. \square A_m) = \Pr(A_1)\Pr(A_2/A_1)...\Pr(A_m/A_1...A_{m-1}) = \Pr(A_1)B_2(x)\Pr(A_2)...B_m(x)\Pr(A_m) = C_m(x)/\log^m x, \qquad (3.16)$$

where $C_m(x) = \prod_{i=2}^{m} B_i(x)$.

Let us denote $\omega_m(p)$ - the number of comparison solutions:

$$x(x+2k_1)...(x+2k_m) \equiv 0 (\bmod\, p).$$

For example, if $m=2$, then $\omega_2(p) = 2$ with $\omega_2(p) = 2$.

Having in mind (3.11) and using this notation, we obtain:

$$C_2(x) = \frac{\Pr(A_2/A_1)}{\Pr(A_2)} \square\, 2\prod_{p \leq \sqrt{x}} \frac{p(p-2)}{(p-1)^2} = \frac{\prod_{p \leq \sqrt{x}} 1-\omega_2(p)/p}{\prod_{p \leq \sqrt{x}} (1-1/p)^2}. \qquad (3.17)$$

Based on (3.16) we get for the value $m=3$:

$$C_3(x) = \frac{\Pr(A_1)\Pr(A_2/A_1)\Pr(A_3/A_1 A_2)}{\Pr(A_1)\Pr(A_2)\Pr(A_3)} = C_2(x)\frac{\Pr(A_3/A_1 A_2)}{\Pr(A_3)}. \qquad (3.18)$$

Having in mind (3.18) we obtain:

$$\frac{C_3(x)}{C_2(x)} = \frac{\Pr(A_3/A_1 A_2)}{\Pr(A_3)}. \qquad (3.19)$$

We continue (3.18) and (3.19) by induction and get:

$$C_m(x) = \frac{\Pr(A_1)\Pr(A_2/A_1)...\Pr(A_m/A_1...A_{m-1})}{\Pr(A_1)\Pr(A_2)...\Pr(A_m)} = C_2(x)\frac{C_3(x)...C_{m-1}\Pr(A_m/A_1...A_{m-1})}{C_2(x)...C_{m-2}\Pr(A_m)}. \quad (3.20)$$

Based on (3.20) and continue (3.17) by induction, we obtain:

$$C_m(x) \approx \frac{0.5e^\gamma \prod_{p \leq x} \frac{1-\omega_m/p}{(1-1/p)^{m-1}}}{0.5e^\gamma \prod_{p \leq x} 1-1/p} = \frac{\prod_{p \leq x} 1-\omega_m/p}{\prod_{p \leq x}(1-1/p)^m}, \quad (3.21)$$

which corresponds to (3.12).

Meaning $C_m = 0$, as starting from some natural number $l$, natural numbers $2(n_1+...+n_l),...,2(n_1+...+n_m)$ are a complete system of deductions.

For example, the numbers $0,2,4$ a complete system of residues modulo 3, so for the triplet $x, x+2, x+4$ the number of comparison solutions $x(x+2)(x+4) \equiv 0 \pmod 3$ is 3 and the value in formula (3.21) is $1 - \frac{\omega_3(p)}{p} = 1 - \frac{3}{3} = 0$, therefore $C_3 = 0$.

Thus, Lemma 3.1 removes the question of the incorrect assumption in the Hardy-Littlewood conjecture of the independence of events that a natural number $x$ is not divisible by primes $p \leq \sqrt{x}$.

The proof of conjectures 3.2 and 3.3 is only their justification, since it is based on assumption (3.3).

The question arises with respect to assumption (3.3). How is such a high accuracy of the conjecture achieved if the probability of a natural number being prime does not exist at all? The following theorem answers the indicated question.

Theorem 3.4

Hardy-Littlewood conjecture is fulfilled, if we replace assumption (3.3) to the existing probability of a natural number from the interval $[2, x)$ to be prime.

Proof

Based on Lemma 2.1, the probability of a natural number from the interval $[2, x)$ to be prime is equal to:

$$P_{2,x}(f) = 1/\log(x)(1+o(1)) = 1/\log(x) + o(1/\log(x)). \qquad (3.22)$$

It can be seen (from a comparison of (3.3) and (3.22)) that the values of the nonexistent and existing probability differ by a value $o(1/\log(x))$.

We denote the sequence of primes k- tuples - $f_k$, and the probability that the k- tuple consists from only primes - $\Pr(f_k)$, then having in mind (3.12) we get:

$$\Pr(f_k) \sim C_k / \log^k(x), \qquad (3.23)$$

where $C_k$ is the constant depending on the type of tuple.

We can write (3.23) as:

$$\Pr(f_k) = C_k(1+o(1))/\log^k(x). \qquad (3.24)$$

If we replace the nonexistent probabilities in (3.24) with the existing ones, then taking into account (3.22), we get:

$$P_{2,x}(f_k) = C_k(1+o(1))(1+o(1))^k / \log^k(x) = C_k(1+o(1))/\log^k(x) = C_k/\log^k(x) + o(1/\log^k(x)), \qquad (3.25)$$

which corresponds to the general form of the Hardy-Littlewood conjecture.

4. DISTRIBUTION OF PRIME TUPLES BASED ON HARDY-LITTLEWOOD CONJECTURE

First, let us consider the arithmetic function of the number of natural numbers that have a certain property - $Q(n)$ such that it can be represented in the following form:

$$Q(n) = \sum_{i=1}^{n} x_i, \qquad (4.1)$$

where a random variable $x_i = 1$ if $i$ is a positive integer which has the indicated property with probability $P(x_i) = P(i)$ and $x_i = 0$ otherwise with probability $P(x_i = 0) = 1 - P(i)$.

Then, based on (4.1), the mathematical expectation $M[x_i] = P(i)$ and, accordingly, the average value of the arithmetic function $Q(n)$ is equal to:

$$M[Q,n] = \sum_{i=1}^{n} M[x_i] = \sum_{i=1}^{n} P(i) \approx \int_{i=1}^{n} P(x)dx. \qquad (4.2)$$

Having in mind that the variance of independent random variables $x_i$ - $D[x_i] = P(i) - P^2(i)$, we obtain that the variance of the arithmetic function $Q(n)$ is equal to:

$$D[Q,n] = \sum_{i=1}^{n} (P(i) - P^2(i)) \approx \int_{i=1}^{n} (P(x) - P^2(x))dx. \qquad (4.3)$$

It was shown in [6] that if independent random variables $x_i$ have finite absolute central moments of the third order and the series $\sum_{i=1}^{n}(P(i) - P^2(i))$ diverges, then the conditions of the central limit theorem in the Lyapunov form are satisfied and $Q(n)$ have an asymptotic normal distribution.

You can use another special case of the central limit theorem [7] to prove that $Q(n)$ has an asymptotic normal distribution.

Let there are independent random variables $x_i$ such that $|x_i| \leq K < \infty$, where $K$ is some constant and $D_n \to \infty, n \to \infty$. There is (in our case) $K = 1, D_n = D[Q,n] = \sum_{i=1}^{n}(P(i) - P^2(i)) \to \infty, n \to \infty$, therefore, the conditions of the special case of the central limit theorem are satisfied. Now let use it for prime tuples.

Denote the density of prime $k$- tuples in the interval of the natural series $[2,n)$ - $A_k(n)$. As I already wrote, there is a probability equal to the density of a strictly increasing sequence on a finite interval of the natural series. Therefore, there is the likelihood of a prime $k$-tuple appearing on the interval of the natural series $[2,n)$:

$$Pr(f_k, 2, n) = A_k(n) = C_k / \log^k(n) + o(1/\log^k(n))$$

or $A_k(n) \sim C_k / \log^k(n)$. \qquad (4.4)

Value $x_i = 1$ if $i$ has the property that there is the first prime number of a prime $k$-tuple.

Having in mind (4.2), (4.3), we can write that the average number of $k$-tuples in the interval of the natural series $[2, n)$ is determined by the asymptotic formula:

$$M_k(n) \sim C_k \int_2^n \log^{-k}(x) dx. \qquad (4.5)$$

Note that (4.5) coincides with the asymptotic formula for the number of $k$-tuples in the Hardy-Littlewood conjecture.

Having in mind that the series $\sum_2^\infty (C_k / \log^k(n) - C_k^2 / \log^{2k}(n))$ diverges, based on the Central Limit Theorem in the Lyapunov form, the number of prime $k$-tuples (under the above assumptions) has an asymptotic normal distribution.

Based on (4.3), we can write, for example, the standard deviation formula for the number of prime twins in the interval of the natural series $[2, n)$:

$$\sigma_2(n) \sim \sqrt{C_2 \int_2^n \log^{-2}(x) dx - C_2^2 \int_2^n \log^{-4}(x) dx}, \qquad (4.6)$$

where $C_2 = 1,32...$

Calculations according to formula (4.6) with the value $n = 10^5$ give $\sigma_2(10^5) = 35$ with the difference between the Hardy-Littlewood conjecture calculated by the formula and the actual number of prime twins $1249 - 1224 = 25$, with the value $n = 10^6$ - $\sigma_2(10^6) = 90$ with a difference $8248 - 8169 = 79$, with the value $n = 10^7$ - $\sigma_2(10^7) = 242$ with a difference $58754 - 58980 = -226$. These calculations fit into the normal distribution.

### 4. CONCLUSION AND SUGGESTIONS FOR FURTHER WORK

The probabilistic estimations can be used in many cases as the only method of analysis of the distribution of primes. Therefore it is opportune to obtain probabilistic estimates of the accuracy of other indicators of the distribution of primes and prime k-tuples.

### 5. ACKNOWLEDGEMENTS

Thanks to everyone who has contributed to the discussion of this paper.